\documentclass[11pt,reqno]{amsart}
\usepackage{amssymb,latexsym}
\usepackage{graphicx}
\usepackage{fancyhdr}

\newcommand{\N}{{\mathbb N}}
\newcommand{\la}{\lambda}

\theoremstyle{plain} \numberwithin{equation}{section}
\newtheorem{thm}{Theorem}[section]
\newtheorem{theorem}[thm]{Theorem}

\renewenvironment{proof}[1][\proofname]{{\it #1.}}{\hfill$\square$}

\begin{document}

\setcounter{page}{1}

\title[ON THE MODES OF THE POISSON DISTRIBUTION OF ORDER K]{On the modes of the Poisson distribution of order k}
\author{Constantinos~Georghiou}
\address{Department of Engineering Sciences\\
                University of Patras\\
                Patras 26500\\
                Greece}
\email{c.georghiou@upatras.gr}
\author{Andreas~N.~Philippou}
\address{Department of Mathematics\\
               University of Patras\\
               Patras 26500\\
               Greece and \vspace{-2.6mm}}
               \address{Technological Educational Institute of Lamia, Lamia, Greece}
\email{anphilip@master.math.upatras.gr}
\author{Abolfazl~Saghafi}
\address{School of Mathematics\\
               Iran University of Science and Technology\\
               Tehran, Iran}
\email{asaghafi@iust.ac.ir}

\begin{abstract}
Sharp upper and lower bounds are established for the modes of the
Poisson distribution of order k. The lower bound established in this
paper is better than the previously established lower bound. In
addition, for $k = 2, 3, 4, 5$, a recent conjecture is presently
proved solving partially an open problem since 1983.
\end{abstract}

\maketitle

\section{Introduction and Summary}
For any given positive integer $k$, denote by $N_k$ the number of
independent trials with constant success probability $p$ until the
occurrence of the $k^{\text{th}}$ consecutive success, and set $q=1
- p$. For $n\geq k$, Philippou and Muwafi \cite{PM} derived the
probability $P(N_k=n)$ in terms of multinomial coefficients and
noted that $P(N_k=n\,|\,p=1/2)=f_{n-k+1}^{(k)}/2^n$ where
$f_n^{(k)}$ is the $n^{\text{th}}$ Fibonacci number of order k. See
also \cite{Feller}, \cite{Shane}, and \cite{Turner}. Philippou et
al. \cite{PGP} showed that $\sum_{n=k}^{\infty} P(N_k=n)=1$ and
named the distribution of $N_k$ \textit{the geometric distribution
of order $k$ with parameter $p$}, since for $k=1$ it reduces to the
geometric distribution with parameter $p$. Assuming that $N_{k,1},
\ldots, N_{k,r}$ are independent random variables distributed as
geometric of order k with parameter $p$, and setting $Y_{k,r}=
\sum_{j=1}^{r} N_{k,j}$, the latter authors showed that
\begin{equation*}\label{Gok}
P(Y_{k,r}=y) = p^{y} \sum
\binom{y_1+\cdots+y_k+r-1}{y_1,\cdots,y_k,r-1}
(\frac{q}{p})^{y_1+\cdots+y_k} \quad y= kr, kr+1, \ldots,
\end{equation*}
where the summation is taken over all $k$-tuples of non-negative
integers $y_1, y_2, \cdots, y_k$ such that $y_1 +2y_2 + \cdots+ ky_k
= y - kr$. They named the distribution of $Y_{k,r}$ \textit{the
negative binomial distribution of order $k$ with parameters $r$ and
$p$}, since for $k = 1$ it reduces to the negative binomial
distribution with the same parameters. Furthermore they showed that,
if $rq\rightarrow \lambda \,(\lambda>0)$ as $r\rightarrow\infty$ and
$q\rightarrow0$, then
\begin{equation}\label{Pok}
\lim_{r\rightarrow\infty}P(Y_{k,r}-kr=x)= \sum_{x_{1},\cdots,x_{k}} \\
e^{-k\lambda} \frac{\lambda^{x_{1}+x_{2}+\cdots+x_{k}}}{x_{1}!\cdots
x_{k}!}= f_{k}(x;\lambda), \quad x=0, 1, 2, \ldots,
\end{equation}
where the summation is taken over all $k$-tuples of non-negative
integers $x_{1}, x_{2}$, $\cdots$, $x_{k}$ such that
$x_{1}+2x_{2}+\cdots+kx_{k}=x$. They named the distribution with
probability mass function $f_k(x;\lambda)$ \textit{the Poisson
distribution of order $k$ with parameter $\lambda$}, since for $k
=1$ it reduces to the Poisson distribution with parameter $\lambda$.
See also \cite{Aki84}, \cite{P1986}, and \cite{Bala}.

Denote by $m_{k,\la}$ the mode(s) of $f_k(x;\la)$, i.e. the value(s)
of $x$ for which  $f_k(x;\la)$ attains its maximum. It is well known
that $m_{1,\la}= \la$ or $\la - 1$ if $\la\in\N$, and
$m_{1,\la}=\lfloor\la\rfloor$ if $\la\not\in\N$. Philippou
\cite{P1983} derived some properties of $f_k(x;\la)$ and posed the
problem of finding its mode(s) for $k \geq2$. See also \cite{SIAM}
and \cite{P2010}.

Hirano et al. \cite{Hirano84} presents several graphs of
$f_{k}(x;\lambda)$ for $\lambda\in(0, 1)$ and $2\leq k\leq8$, and
Luo \cite{Luo} derived the following inequality
\begin{equation}
m_{k,\lambda}\geq k\lambda \ ^{k}\sqrt{k!}-\frac{k(k+1)}{2},\quad
k\geq1 \ (\lambda>0),
\end{equation}
which is sharp in the sense that $m_{1,\lambda}=\lambda-1$ for
$\la\in\N$. Recently, Philippou and Saghafi \cite{SSIAM} conjectured
that, for $k\geq2$ and $\lambda\in\mathbb{N}$,
\begin{equation}\label{PhiSa}
m_{k,\lambda}= \lambda k(k+1)/2- \lfloor k/2\rfloor,
\end{equation}
where $\lfloor u\rfloor$ denotes the greatest integer not exceeding
$u\in\mathbb{R}$.

In this paper, we employ the probability generating function of the
Poisson distribution of order $k$ to improve the bound of Luo
\cite{Luo} and to also give an upper bound (see Theorem 2.1). We
then use Theorem 2.1 to prove the conjecture of Philippou and
Saghafi \cite{SSIAM} when $k = 2, 3, 4, 5$, partially answering the
open problem of Philippou \cite{P1983,SIAM,P2010}.

\section{Main results}
In the present section, we state and prove the following two
theorems.
\begin{theorem}\label{th:2l}
For any integer $k\geq1$ and real $\lambda>0$, the mode(s) of the
Poisson distribution of order $k$ satisfies(y) the inequalities
\begin{equation*}
\lfloor \lambda k(k+1)/2\rfloor - \frac{k(k+1)}{2} +1- \delta_{k,1}
\leq m_{k,\lambda}\leq \lfloor\lambda k(k+1)/2\rfloor,
\end{equation*} where $\delta_{k,1}$ is the Kronecker delta.
\end{theorem}

\begin{theorem}\label{th:22}
For $\lambda\in\mathbb{N}$ and $2\leq k \leq5$, the Poisson
distribution of order k has a unique mode $m_{k,\lambda}=\lambda
k(k+1)/2 -\lfloor k/2 \rfloor$.
\end{theorem}

For the proofs of the theorems we employ the probability generating
function of the Poisson distribution of order k and some recurrences
derived from it. We observe first that the left hand side inequality
in Theorem 2.1 is sharp since, for $\la\in\N, m_{1,\la}= \la-1$, the
value of the lower bound for $k = 1$. The right hand side inequality
is also sharp in the sense that there exist values of $k$ and $\la$
for which $m_{k,\la}= \lfloor\la k(k+1)/2\rfloor$. We also note that
our lower bound is better than that of Luo \cite{Luo} for $k\geq 2$.
\vspace{2.5mm}

\noindent\begin{proof}[Proof of Theorem 2.1] For notational
simplicity, we presently set $P_{x}= f_{k}(x;\lambda)$, omitting the
dependence on $k$ and $\lambda$, and $\Delta_x = P_x - P_{x-1}$, $x
= 0, 1, \ldots$. It is easily seen \cite{Galkin,P1983,P1988} that
the probability generating function of $P_x$ is
\begin{equation}\label{gn}
g(s)=
\sum_{x=0}^{\infty}s^{x}P_{x}=e^{\lambda(-k+s+s^{2}+\cdots+s^{k})},
\end{equation}
which implies that
\begin{equation}\label{gpn}
g^{\prime}(s)= \lambda(1+2s+\cdots+ks^{k-1}) g(s).
\end{equation}

For $x \geq1$, we differentiate $(x-1)$ times $g^{\prime}(s)$ and
employ the fact that $P_x=(\frac{1}{x!})\frac{\partial^x
g(s)}{\partial s^x}$ at $s = 0$ to get the recurrence
\begin{equation}\label{xPx}
xP_{x}= \sum_{j=1}^{k} j\lambda P_{x-j}, \quad x \geq1.
\end{equation}

We note that (\ref{xPx}) is trivially true for $x = 0$.  By
definition $P_x\leq P_{m_{k,\la}}$ for every $x \geq 0$, and
therefore
\begin{equation*}
xP_{x}=\sum_{j=1}^{k} j\lambda P_{x-j} \leq \sum_{j=1}^{k} j\lambda
P_{m_{k,\lambda}}= \lambda P_{m_{k,\lambda}}k(k+1)/2.
\end{equation*}
Upon setting $x = m_{k,\lambda}$ we get $m_{k,\la} \leq \la
k(k+1)/2$. Therefore $m_{k,\la}\leq \lfloor\lambda k(k+1)/2\rfloor$
since $m_{k,\la}$ is a non negative integer.

As for the left hand side inequality we note that it is trivially
true for $k = 1$ and $\la > 0$, since $m_{1,\la}= \la$ or $\la - 1$
if $\la\in\N$, and $m_{1,\la}=\lfloor\la\rfloor$ if $\la\not\in\N$.
Therefore we assume that $k \geq 2$.  For $0 < \la < 1$, the
inequality is true since $\lfloor \lambda k(k+1)/2\rfloor -
\frac{k(k+1)}{2} +1 \leq0$. For $\la = 1$ it is also true since
$e^{-k}=P_0 = P_1 < P_2 = 3e^{-k}/2$. Let then $\la > 1$. We will
show that $P_x$  increases, or, equivalently, $\Delta_x$ is
positive, for $0 \leq x \leq \lfloor \la
k(k+1)/2\rfloor-\frac{k(k+1)}{2}+1$.

From the definition of $\Delta_x$ and (\ref{gn}), we obtain
\begin{equation}\label{dn}
h(s)= \sum_{x=0}^{\infty}s^{x}\Delta_{x}=(1-s)g(s).
\end{equation}
Differentiating $h(s)$ twice we get
\begin{eqnarray}\begin{array}{ll}\label{hpps}
h^{\prime\prime}(s)=
\lambda\left(\lambda\sum_{j=1}^{k-1}\frac{j(j+1)}{2}s^{j-1}
+\left(\frac{k(k+1)}{2}\right)(\lambda-2)s^{k-1}+s^{k}f(s)\right)g(s),
\end{array}\end{eqnarray}
where $f(s)=\sum_{j=0}^{k-1}a_{j}s^{j}$ is a $(k-1)^{\text{th}}$
degree polynomial. Next, differentiating $x$ times
$h^{\prime\prime}(s)$ and then setting $s=0$, we get
\begin{equation*}\begin{array}{ll}
\frac{(x+1)(x+2)}{\lambda}\Delta_{x+2}=\sum_{j=1}^{k-1}\frac{j(j+1)}{2}\lambda
P_{x+1-j}+\frac{k(k+1)}{2}(\lambda-2)P_{x-k+1}+\sum_{j=0}^{k-1}a_{j}P_{x-k-j}.
\end{array}\end{equation*}
Finally, eliminating successively
$P_{x-2k+1},P_{x-2k},\cdots,P_{x-k}$, by means of (\ref{xPx}) we
arrive at
\begin{equation}\begin{array}{ll}\label{Dx2}
\frac{(x+1)(x+2)}{\lambda}\Delta_{x+2}=
\sum_{j=1}^{k-1}\left(j\lambda+x+1-j\right)P_{x+1-j}
+k(\lambda-x-2)P_{x-k+1}.
\end{array}\end{equation}

Setting $P_{x+1-j}=\Delta_{x+1-j}+P_{x-j}$ in (\ref{Dx2}) we obtain

\begin{equation}\begin{array}{l}\label{Dx1}\vspace{1.5mm}
\frac{(x+1)(x+2)}{\lambda}\Delta_{x+2}= \sum_{j=1}^{k-1}\left(j(x+1)
+\frac{(\lambda-1)j(j+1)}{2}\right) \Delta_{x+1-j} +
\left(\frac{(\lambda-1)k(k+1)}{2}-1-x\right) P_{x+1-k}.
\end{array}\end{equation}

Since $\la>1$, we have $\Delta_0 = e^{-k\la} > 0$ and $\Delta_1 =
(\la-1) e^{-k\la} >0$. An easy recursion using (\ref{Dx1}) shows
that $\Delta_x>0$ for $2\leq x+2 \leq\frac{(\lambda-1)k(k+1)}{2}+1$
also. This completes the proof of the theorem.
\end{proof} \vspace{2.5mm}

\noindent\begin{proof}[Proof of Theorem 2.2] For $k = 2$, Theorem
2.1 reduces to $3\la-2 \leq m_{2,\la} \leq 3\la$. Therefore, in
order to show that $m_{2,\la}=3\la-1$, it suffices to show that
$\Delta_{3\la-1}>0$ and $\Delta_{3\la}<0$. However, by (\ref{xPx}),
$3\Delta_{3\la}=-2\Delta_{3\la-1}$. Therefore, we will only show
$\Delta_{3\la-1} > 0$. For $\lambda = 1$, $\Delta_{3\la-1}=\Delta_2=
e^{-2}/2 > 0$; for $\la = 2$, $\Delta_{3\la-1}= \Delta_5 =
\frac{4e^{-4}}{15}>0$. Let $\la\geq3$ and $x = 3\la - 3$. Using
(\ref{Dx2}) we have
\begin{equation*}\begin{array}{ll}
\frac{(3\lambda-1)(3\lambda-2)}{\lambda} \Delta_{3\lambda-1}=
(4\lambda-3) P_{3\lambda-3}-(4\lambda-2)P_{3\lambda-4}=
(4\la-3)\Delta_{3\la-3}-P_{3\la-4}.
\end{array}\end{equation*}
By (\ref{xPx}),
\begin{equation*}\begin{array}{ll}
\frac{1}{\la^3}\prod_{j=1}^{6}(6\la-j)\Delta_{3\la-1}=
(64\la^3-267\la^2+360\la-156) \Delta_{3\la-7}+3(\la^2+8\la-12)
P_{3\la-8}.
\end{array}\end{equation*}
Therefore, $\Delta_{3\la-1}$ is positive since $\Delta_{3\la-7}>0$
by Theorem 2.1, $P_{3\la-8}>0$ by (1.1), and both
$64\la^3-267\la^2+360\la-156$ and $\la^2+8\la-12$ take positive
values.\\

For $k = 3$, Theorem 2.1 reduces to $6\la-5 \leq m_{3,\la} \leq
6\la$. Therefore, in order to show that $m_{3,\la}=6\la-1$, it
suffices to show that $\Delta_{6\la-j}>0 \ (1\leq j \leq4)$ and
$\Delta_{6\la}<0$. However,
$6\Delta_{6\la}=-5\Delta_{6\la-1}-3\Delta_{6\la-2}$ because of
(\ref{xPx}). We will show then only that $\Delta_{6\la-4}>0$ (the
other three can be treated similarly). For $\la = 1$,
$\Delta_{6\la-4}= \Delta_2 = e^{-3}/2>0$. Let $\la\geq2$ and $x =
6\la - 6$. Using (\ref{Dx2}) we have

\begin{equation*}\begin{split}
\frac{(6\lambda-4)(6\lambda-5)}{\lambda} \Delta_{6\lambda-4}&=
(7\lambda-6) P_{6\lambda-6}+(8\lambda-7)P_{6\lambda-7}-
(15\lambda-12)P_{6\lambda-8}\\
&= (7\la-6) \Delta_{6\la-6}+(15\la-13) \Delta_{6\la-7}-P_{6\la-8}.
\end{split}\end{equation*}
By (\ref{xPx}),
\begin{equation*}\begin{array}{ll}\vspace{2mm}
\frac{1}{\la^3}&\!\!\!\prod_{j=4}^{8}(6\la-j)\Delta_{6\la-4}=
(1015\la^3-3234\la^2+3396\la-1176) \Delta_{6\la-9}\\
&+(1203\la^3-3610\la^2+3576\la-1176)\Delta_{6\la-10}
+2(199\la^2-372\la+168)P_{6\la-11},
\end{array}\end{equation*}
which is positive, since $\Delta_{6\la-9}>0$ and
$\Delta_{6\la-10}>0$ by Theorem 2.1, $P_{6\la-11}>0$ by (1.1),  and
their
polynomial coefficients take positive values. \\

When $k = 4 \ (k = 5)$ we use the same procedure as above to show
that $\Delta_{10\la-j}>0 \ (2 \leq j \leq 8)$ and
$\Delta_{10\la-1}<0 \ (\Delta_{15\la-j}>0 \ (2 \leq j \leq 13)$ and
$\Delta_{15\la-1}<0)$. Therefore, $m_{4,\la}=10\la-2 \
(m_{5,\la}=15\la-2)$.\\

{\bf Remark 2.1.} As $k$ increases the computations become
increasingly difficult and lengthy. We have used the computer
algebra system Derive and a personal computer to check them.\\

{\bf Remark 2.2.} According to the conjecture of Philippou and
Saghafi \cite{SSIAM}, $m_{6,2} = 39$. However, by (\ref{xPx}) (and
(\ref{Pok})), we presently find that $f_6(40;2)=0.0297464817
>0.0297385179=f_6 (39;2)$. Therefore the conjecture is not true at
least for $k = 6$ and $\la=2$.
\end{proof}

\section{Further research}
In this note we have derived an upper and a lower bound for the
mode(s) of the Poisson distribution of order $k$. Our lower bound is
better than that of Luo \cite{Luo}. We have also established the
conjecture of Philippou and Saghafi \cite{SSIAM} for $2 \leq k \leq
5$ and $\la\in\N$, partially solving the open problem of Philippou
\cite{P1983,SIAM,P2010}. However, the problem remains open for all
other cases.\vspace{5mm}

\noindent{\bf Acknowledgement.} The authors would like to thank the
anonymous referee for very helpful comments on the style of this
paper.

\medskip

\noindent MSC2010: 60E05, 11B37, 39B05

\end{document}